\documentclass[11pt]{article}

\usepackage{tikz}
\usepackage{subfigure}
\usepackage[english]{babel}

\usepackage[center]{caption2}
\usepackage{amsfonts,amssymb,amsmath,latexsym,amsthm}
\usepackage{multirow}

\newtheorem{thm}{Theorem}[section]
\newtheorem{cor}[thm]{Corollary}
\newtheorem{lem}[thm]{Lemma}
\newtheorem{prop}[thm]{Proposition}

\newtheorem{quest}[thm]{Question}

\newtheorem{definition}{Definition}

\begin{document}
	
\title{On Induced Subgraph of Cartesian Product of Paths}
\author{Jiasheng Zeng$^a$, \quad Xinmin Hou$^{b,c}$\footnote{Email address: xmhou@ustc.edu.cn (X. Hou)}\\
\small $^{a}$ School of the Gifted Young\\
\small University of Science and Technology of China, Hefei, Anhui 230026, China.\\
\small $^{b}$ School of Mathematical Sciences\\
\small University of Science and Technology of China, Hefei, Anhui 230026, China.\\
\small  $^{c}$ CAS Key Laboratory of Wu Wen-Tsun Mathematics\\
\small University of Science and Technology of China, Hefei, Anhui 230026, China.}
	
	\date{}
	
	\maketitle

\begin{abstract}
Chung, F\"uredi, Graham, and Seymour (JCTA, 1988) constructed an induced subgraph of the hypercube $Q^n$ with $\alpha(Q^n)+1$ vertices and with maximum degree smaller than $\lceil \sqrt{n} \rceil$. Subsequently, Huang (Annals of Mathematics, 2019) proved the Sensitivity Conjecture by demonstrating that the maximum degree of such an induced subgraph of hypercube $Q^n$ is at least $\lceil \sqrt{n} \rceil$, and posed the question: Given a graph $G$, let $f(G)$ be the minimum of the maximum degree of an induced subgraph of $G$ on $\alpha(G)+1$ vertices,  what can we say about $f(G)$?  In this paper, we investigate this question for Cartesian product of paths $P_m$, denoted by $P_m^k$. We determine the exact values of $f(P_{m}^k)$ when $m=2n+1$ by showing that $f(P_{2n+1}^k)=1$ for $n\geq 2$ and $f(P_3^k)=2$, and give a nontrivial lower bound of $f(P_{m}^k)$ when $m=2n$ by showing that $f(P_{2n}^k)\geq \lceil \sqrt{\beta_nk}\rceil$. In particular, when $n=1$, we have $f(Q^k)=f(P_{2}^k)\ge \sqrt{k}$, which is  Huang's result. The lower bounds of $f(P_{3}^k)$ and $f(P_{2n}^k)$ are given by using the spectral method provided by Huang.
\end{abstract}

\section{Introduction}

We consider only simple and finite graph in this paper. For a graph $G=(V,E)$ and a vertex $v\in V$, write $d_G(v)$ for the degree of $v$, the maximum degree of $G$  is $\Delta(G)=\max\{d_G(v) : v\in V\}$. For a subset $S\subseteq V$, write $G[S]$ for the subgraph induced by $S$.  Let $\alpha(G)$ be the independence number of $G$, i.e. the maximum size of an independent set of $G$. Define
$$f(G)=\min\{ \Delta(G[S]) : S\subset V(G) \text{ with } |S|=\alpha(G)+1\}.$$ 
Let $G$ and $H$ be two graphs.  The {\em Cartesian product} of graphs $G$ and $H$, denoted by $G\square H$, is the graph with vertex set $V(G) \times V(H) $, in which two vertices, say $(x_1,x_2)$ and $(y_1,y_2)$, are adjacent if and only if $x_1=y_1$ and $x_2$ is adjacent to $y_2$ in $H$, or $x_2=y_2$ and  $x_1$ is adjacent to $y_1$ in $G$. Write $G^k=G\square G\square \cdots \square G$ for the  Cartesian product of $k$ copies of graph $G$. For example, the well known hypercube $Q^n=K_2^n=K_2\square K_2\square\cdots\square K_2$.

A signed graph $(\Gamma, \sigma)$ of a graph $G(V,E)$ is a graph with vertex set $V$ and edge set $E$, together with a map $\sigma:E\longrightarrow\{-1,+1\}$. The {\em adjacent matrix} of $\Gamma$ is defined to be a symmetric $\{0,\pm1\}$ matrix $A(\Gamma)=(a_{ij}^{\sigma})$ with $a_{ij}=\sigma(v_iv_j)$.
Clearly, $A(\Gamma)$ is a matrix by putting some minus on the adjacent matrix of $G$, we call it a signed matrix of $G$. For any matrix $A$, define $\lambda_i(A)$ as the $i$-th largest eigenvalue of $A$.

 Chung et al.~\cite{Chung} provided a construction proof to show that $f(Q^n)\leq \sqrt{n}$ in 1988, and Huang~\cite{Huang} proved that (Huang-Theorem) $f(Q^n)\geq \lceil \sqrt{n} \rceil$  by defining a group of signed graphs of $Q^n$. Huang-Theorem is well known since then as the fact that it is equivalent to the Sensitivity Conjecture (for readers interested, one can refer to two excellent surveys provided by Hatami, Kulkarni, and Pankratov~\cite{AboutConjectur} and  by Karthikeyan, Sinha, and Patil~\cite{S.Sinha}, respectively).
In~\cite{Huang},  Huang also proposed the following interesting problem.
\begin{quest}[Huang, \cite{Huang}]\label{QUES: q1}
	Given a ``nice" graph $G$ with high symmetry,   what can we say about $f(G)$? In particular,
	for which graphs, the method used in proving Huang-Theorem would provide a tight
	bound?
\end{quest}

There are some results in this flavor and some extensions of Huang's Theorem, for example, in~\cite{Alon,Tsang,Verret} for Cayley graphs, Tikaradze~\cite{Tikara} for the Cartesian product of oriented cycles, Hong, Lai and Liu~\cite{CartesianType} for product of signed bipartite graphs.
The Hamming graph $H(n,k)=K_k^n=K_k\square K_k\ldots\square K_k$ is a natural extension of the hypercube $Q^n=K_2^n$, where $K_k$ is a complete graph on $k$ vertices.  Recently, Dong~\cite{Dong} extended the result of Chung et al~\cite{Chung} to Hamming graph by showing that  $f(H(n,k))\leq \lceil \sqrt{n} \rceil$ for $k\geq 3$. More recently, Tandya~\cite{Vincent} improved Dong's result by showing that $f(H(n,k))=1$ for all $k\geq 3$. 


In this paper, we consider Question~\ref{QUES: q1} when $G=P_m^k$, the Cartesian product of paths $P_m$, where $P_m$ is a path on $m$ vertices. 
The main result is the following. 


\begin{thm}\label{main1}
Let $k\ge 1$ be an integer. Then
$$
f(P_m^k)\begin{cases}
	=2, & m=3\\
	=1, &  m=2n+1\geq 5\\
	\geq \lceil \sqrt{\beta_nk}\rceil, & m=2n\ge 2
	\end{cases},
$$
where $\beta_n>0$ is a constant depending on $n$. In particular, $\beta_1 = 1$.
\end{thm}

Note that when $n=1$, $P_{2n}^k=P_2^k=Q^k$ and $\beta_1=1$. Therefore, Huang-Theorem is a direct corollary of Theorem~\ref{main1}.
\begin{cor}[Huang~\cite{Huang}]\label{cor1}
$$f(P_2^n)=f(Q^n)\geq \lceil \sqrt{n} \rceil.$$ 
\end{cor}

The rest of the article is arranged as follows. In Section 2, we give some preliminaries and lemmas. We prove Theorem~\ref{main1} in Section 3 when $m$ is odd and in Section 4 when $m$ is even. We give some discussion and Remarks in the last section. 



\section{Preliminaries and lemmas}

It will be convenient to view the vertices of path $P_n$ as $[n]:=\{1,2,\cdots,n\}$. So the vertices of $P_n^k$ can be written as $[n]^k:=[n]\times\cdots\times[n]$. 
Write $u\sim v$ if $u$ and $v$ are adjacent in $P_n^k$ and $u\nsim v$, otherwise. 

\begin{definition} 
	For  $a\in [n]$, define a map $Q_{n,k+1}^{a}:\mbox{ }[n]^k \longrightarrow [n]^{k+1}$ with 	$$Q_{n,k+1}^{a}((x_1,\cdots,x_k))=(x_1,\cdots,x_k,a). $$
\end{definition}

The following is a simple observation from the definition of  $P_{n}^k$. 
\begin{prop}\label{neighborcon}
	Let $u=(u_1,\ldots, u_k), v=(v_1,\ldots, v_k)\in [n]^k$ be two vertices of $P_{n}^k$. Then $u\sim v$ in $P_{n}^k$ if and only if the 1-norm  $||u-v||_1=\sum_{i=1}^k|u_i-v_i|=1$.
\end{prop}

The following lemma shows that the map $Q_{n,k+1}^a$ preserve the independence of its preimage.
\begin{lem}\label{propQ}
If $I\subseteq [n]^k$ is an independent set of $P_n^k$, then $Q_{n,k+1}^a(I)$ is an independent set of $P_n^{k+1}$ for all $a\in [n]$.
\end{lem}
\begin{proof}
Let $u=(u_1,\cdots,u_k,a)$ and $v=(v_1,\cdots,v_k,a)$ be two vertices in $Q_{n,k+1}^a(I)$. Then $u'=(u_1,\cdots,u_k)$, $v'=(v_1,\cdots,v_k)\in I$, i.e., $u' \nsim v'$ in $P_n^k$. 
By Proposition~\ref{neighborcon}, $\sum_{i=1}^k|u_i-v_i|>1$. Hence, $\sum_{i=1}^k|u_i-v_i|+|a-a|>1$. Therefore, $(u_1,\cdots,u_k,a)\nsim (v_1,\cdots,v_k,a)$ in $P_{n}^{k+1}$. 
\end{proof}

Furthermore, we are likely to use the following theorems to build a bridge between the maximum degree of graphs and the eigenvalues of a signed graph.
\begin{lem}[Cauchy Interlace Theorem~\cite{Cauchy}]\label{Cauchy Interlace Thm}
	Let $A$ be a symmetric $n\times n$ matrix and $B$ be a $m\times m$ principle submatrix of $A$. If the eigenvalues of $A$ are $ {\lambda}_1\geq\cdots\geq{\lambda}_n$ and the eigenvalues of $B$ are $\mu_1\geq\cdots\geq\mu_m$, then for all $1\leq i\leq m$,$$\lambda_i\geq\mu_i\geq\lambda_{n-m+i}.$$
\end{lem}

\begin{lem}[Huang~\cite{Huang}]\label{Huang1}
	Suppose $H$ is an $m$-vertex undirected graph, and $A$ is a symmetric matrix whose entries are in $\{0,\pm1\}$ and whose rows and columns are indexed by $V(H)$, and whenever $u$ and $v$ are non-adjacent in $H$, $A_{uv}=0.$ Then
	$$\Delta(H)\geq\lambda_1(A).$$ 
\end{lem}

Some properties about Cartesian product are also needed:
\begin{lem}[Germina et al.~\cite{Eigenvalues}]\label{product}
	Let A be an $n\times n$ matrix and $B$ be an $m\times m$ matrix. If the eigenvalues of A are $\lambda_1,\ldots,\lambda_n$ and the eigenvalues of $B$ are $\mu_1,\ldots,\mu_m$, then all eigenvalues of $I_m\otimes A+B\otimes I_n$ are $\{\lambda_i+\mu_j,1\leq i\leq n, 1\leq j\leq m\}$.
\end{lem}

\begin{lem}[Hammack et al.~\cite{Bipartite}]\label{Bipartite1}
The Cartesian product graph	$G_1\square G_2$ is bipartite if and only if $G_1$ and $G_2$ are bipartite.
\end{lem}

As a direct corollary of Lemma~\ref{Bipartite1}, we have the following proposition. 
\begin{prop}\label{CPGBi}
	$P_n^k$ is a bipartite graph. 
\end{prop}

The following lemma shows that the non-zero eigenvalues of signed matrix of a bipartite graph are pairwise symmetry with respect to $0$.
\begin{lem}\label{bipartite}
	Suppose that $(\Gamma,\sigma)$ is a signed graph of a bipartite graph $G$.
	Then all of the non-zero eigenvalues of $A(\Gamma)$ are pairwise symmetry with respect to $0$.
\end{lem}
\begin{proof}
	According to the definition of a bipartite graph and the adjacent matrix of a signed graph, $A=A(\Gamma)$ is symmetric and can be presented as $A=
	\begin{pmatrix}
		O & B \\
		B^T & O
	\end{pmatrix} $ if we relabel the vertex of $G$. Suppose $A$ has an eigenvalue $\lambda\neq 0$ and $x=(x_1,x_2)^T\in R^{|V(G)|}$ is one of its eigenvector. Then
	$$Ax=\lambda x, \mbox{ i.e. } 
	\begin{cases}
		Bx_2=\lambda x_1 \\
		B^Tx_1 = \lambda x_2
	\end{cases}.
	$$
	Now let $x'=(x_1,-x_2)^T$. Then
	$$
	Ax'= 
	\begin{pmatrix}
		-Bx_2 \\
		B^Tx_1
	\end{pmatrix}
	=
	\begin{pmatrix}
		-\lambda x_1\\
		\lambda x_2
	\end{pmatrix}
	=-\lambda x',
	$$
i.e.,  $-\lambda$ is also a non-zero eigenvalue of $A$.
\end{proof}

Now, we determine the independence number $\alpha(P_{m}^k)$. 
\begin{definition}
{Let $V_1=\{i\in [m]: i\equiv 1\pmod2\}\subset V(P_{m})$. For $k\ge 2$, we recursively define 
	$$V_{k}=Q_{m,k}^1(V_{k-1})\cup Q_{m,k}^2(\overline{V_{k-1}})\cup \cdots \cup Q_{m,k}^{m}(V_{k-1}), \text{ if $2\nmid m$}$$ or 	$$V_{k}=Q_{m,k}^1(V_{k-1})\cup Q_{m,k}^2(\overline{V_{k-1}})\cup \cdots \cup Q_{m,k}^{m}(\overline{V_{k-1}}), \text{ if $2\mid m$}$$
			where $\overline{V_{k-1}}=[m]^{k-1}-V_{k-1}$.}

\end{definition}

In fact, the vertex set $V_k$ defined above is an independent set of $P_{m}^k$.
\begin{lem}\label{Vk}
The vertex sets	$V_k$ and $\overline{V_k}$ are independent sets of $P_{m}^k$, respectively. Moreover, $\left|V_k \right|=\left\lceil\frac{m^k}{2}\right\rceil$ and $\left|\overline{V_k}\right|=\left\lfloor\frac{m^k}{2}\right\rfloor$.
\end{lem}
\begin{proof}
We just give the proof when $m=2n+1$, the other case can be proved similarly. 
By induction on $k$. When $k=1$, $V_1=\{1,3,\ldots,2n+1\}$ and $\overline{V_1}=\{2,4,\ldots,2n\}$ are independent sets of $P_{2n+1}$. The result holds trivially. Now assume that $V_{k-1}$ and $\overline{V_{k-1}}$ are independent sets of $P_{2n+1}^{k-1}$ for $k\geq 2$ and $\left|V_{k-1} \right|=\frac{(2n+1)^{k-1}+1}{2}$. Then $|\overline{V_{k-1}}|=\frac{(2n+1)^{k-1}-1}{2}$. Therefore, 
$$\left|V_k\right|=(n+1)\cdot \left|V_{k-1}\right|+n\cdot \left|\overline{V_{k-1}}\right|=\frac{(2n+1)^k+1}{2}.$$
According to Lemma~\ref{propQ}, $Q_{2n+1,k}^i(V_{k-1})$ and $Q_{2n+1,k}^i(\overline{V_{k-1}})$ are both independent sets of $P_{2n+1}^k$ for $i\in[2n+1]$. Furthermore, for all $1\leq i\neq j\leq n+1$, there is no edge between $Q_{2n+1,k}^{2i-1}(V_{k-1})$ and $Q_{2n+1,k}^{2j-1}(V_{k-1})$ since $||u-v||_1\ge 2$
for any $u\in Q_{2n+1,k}^{2i-1}(V_{k-1})$ and $v\in Q_{2n+1,k}^{2j-1}(V_{k-1})$. With similar reason there is no edge between $Q_{2n+1,k}^{2i}(\overline{V_{k-1}})$ and $Q_{2n+1,k}^{2j}(\overline{V_{k-1}})$, and between $Q_{2n+1,k}^{2i-1}({V_{k-1}})$ and $Q_{2n+1,k}^{2i}(\overline{V_{k-1}})$.
Therefore, $V_k$ is an independent set of $P_{2n+1}^k$.

With the same discussion, we have 
$$\overline{V_k}=Q_{2n+1,k}^1(\overline{V_{k-1}})\cup Q_{2n+1,k}^2(V_{k-1})\cup \cdots \cup Q_{2n+1,k}^{2n+1}(\overline{V_{k-1}})$$
is also an independent set of $P_{2n+1}^k$ and $$\left|\overline{V_k}\right|=(n+1)\cdot \left|\overline{V_{k-1}}\right|+n\cdot \left|{V_{k-1}}\right|=\frac{(2n+1)^k-1}{2}.$$ 
\end{proof}

The following theorem shows that $V_k$, in fact, is a maximum independent set of $P_{m}^k$.
\begin{thm}\label{OddIn}
	The independence number of $P_{m}^k$ is $\lceil\frac{m^k}{2}\rceil$.
\end{thm}
\begin{proof}
By Lemma~\ref{Vk}, we have $\alpha(P_{m}^k)\geq \left| V_k \right|=\lceil\frac{m^k}{2}\rceil$. Now we show that  $\alpha(P_{m}^k)\leq \lceil\frac{m^k}{2}\rceil$ by introduction on $k$. When $k=1$, we are done from the fact that $\alpha(P_{m})=\lceil\frac{m}{2}\rceil$. Assume that all vertex sets of $V(P_{m}^s)$ with more than $\lceil\frac{m^s}{2}\rceil$ vertices have two adjacent vertices for $1\leq s\leq k-1$.

Now suppose $V\subseteq V(P_{m}^k)$ and $\left|V\right|\geq \lceil\frac{m^k}{2}\rceil+1.$ Define
$$
V(a) = \{(1,a),(2,a),\cdots,(m,a)\} \mbox{, }a\in [m]^{k-1}.
$$
Then $V(a)$ induces a path with $m$ vertices in $P_{m}^k$. We write $V(a)$ for $P_{m}^k[V(a)]$ for convenience in the following. Note that 
$$V(P_{m}^k)=\bigcup_{a\in [m]^{k-1}}V(a)\mbox{ and  }V(a)\cap V(b)=\emptyset\mbox{ for all distinct } a, b\in[m]^{k-1}.$$ 
If $\left|V\cap V(a)\right|\leq \lfloor\frac{m}{2}\rfloor,\mbox{ for all } a\in [m]^{k-1}$, then  
$$
\begin{aligned}
	\left|V \right|&=\sum_{a\in [m]^{k-1}}\left|V\cap V(a)\right| \leq (m)^{k-1}\cdot \left\lfloor\frac{m}{2}\right\rfloor<\left\lceil\frac{m^k}{2}\right\rceil+1,
\end{aligned}
$$
a contradiction.
Thus there must exist $a\in[m]^{k-1}$ with $|V\cap V(a)|\ge \lfloor\frac{m}{2}\rfloor+1$.

For $m=2n$, such a $V\cap V(a)$ with $|V\cap V(a)|\ge n+1$ cannot be an independent set because $\alpha(V(a))=n$. Thus $V$ is not an independent set of $P_{2n}^k$ and we have $\alpha(P_{2n}^k)=\lceil\frac{m^k}2\rceil=n(2n)^{k-1}$. 

Now assume $m=2n+1$. If there exists an $a\in[2n+1]^{k-1}$ such that $\left|V\cap V(a)\right|\geq n+2$, then there are two adjacent vertices in $V\cap V(a)$ because $\alpha(V(a))=n+1$. Thus there must exist $a\in [2n+1]^{k-1}$ with  $\left|V\cap V(a)\right|=n+1$ and $V\cap V(a)$ is an independent set of $V(a)$, we assume that there are $\ell$ such $a\in [2n+1]^{k-1}$, say $\{a_1, a_2,\ldots, a_\ell\}$, i.e. $V\cap V(a_t)=\{(1,a_t),(3,a_t),\cdots,(2n+1,a_t)\}$ for all $t\in[\ell]$.   
Since 
$$
\begin{aligned}
\frac{(2n+1)^k+1}{2}+1&\le\left|V\right| 
	=\sum_{a\in [2n+1]^{k-1}}\left|V\cap V(a)\right| \\
	&\leq \ell\cdot (n+1)+((2n+1)^{k-1}-\ell)\cdot n\\
	&=n(2n+1)^{k-1}+\ell,
\end{aligned}
$$
we have $\ell\geq \frac{(2n+1)^{k-1}+1}{2}+1.$
Define$$\tilde{V}=\{(1,a_1),(1,a_2),\cdots,(1,a_\ell)\}.$$ 
Then $\tilde{V}\subseteq V$.
Let $U=\{(1,a) : a\in [2n+1]^{k-1}\}$. Then $G=P_{2n+1}^k[U]$ is an induced subgraph of $P_{2n+1}^k$ and $G\cong P_{2n+1}^{k-1}$. Note that $|\tilde{V}|=\ell\geq \frac{(2n+1)^{k-1}+1}{2}+1$. By inductive hypothesis,  there are two vertices, say $u=(1,a_i)$ and $v=(1,a_j)$, such that $u$ and $v$ are adjacent in $G\subseteq P_{2n+1}^k$. Therefore, we have $\alpha(P_{2n+1}^k)=|V|-1=\lceil\frac{{(2n+1)}^k}{2}\rceil$.


\end{proof}

\section{Proof of Theorem~\ref{main1} when $m=2n+1$}

Inspired by the construction of $V_k$ in $P_{2n+1}^k$, we define a vertex set in $P_{2n+1}^k$ recursively with size $\alpha(P_{2n+1}^k)+1$.
\begin{definition}
Let
$$
X_1=\{2,4,\cdots,2n\}\cup\{1,2n+1\}\subseteq V(P_{2n+1}),
$$
and, for $k\ge 2$, define
$$
X_k=Q_{2n+1,k}^1(X_{k-1})\cup Q_{2n+1,k}^2(\overline{X_{k-1}})\cup\cdots\cup Q_{2n+1,k}^{2n+1}(X_{k-1}).
$$
\end{definition}

We first show that $\left|X_k\right|=\alpha(P_{2n+1}^k)+1=\frac{(2n+1)^k+1}{2}+1$.
\begin{prop}\label{properX1}
$\left|X_k\right|=\alpha(P_{2n+1}^k)+1=\frac{(2n+1)^k+1}{2}+1.$
\end{prop}
\begin{proof}
We prove by induction on $k$. When $k=1$, we have $\left|X_1\right|=n+2=\alpha(P_{2n+1}^1)+1$. Now suppose that for $k\geq 2$, $\left|X_{k-1}\right|=\alpha(P_{2n+1}^{k-1})+1=\frac{(2n+1)^{k-1}+1}{2}+1$. Then
$$
\begin{aligned}
	\left|X_k\right|
	&=(n+1)\cdot \left|X_{k-1}\right|+n\cdot \left|\overline{X_{k-1}}\right|\\
	&=n\cdot (2n+1)^{k-1}+\left|X_{k-1}\right|\\
	&=n\cdot (2n+1)^{k-1}+\frac{(2n+1)^{k-1}+1}{2}+1\\
	&=\frac{(2n+1)^k+1}{2}+1.
\end{aligned}
$$
\end{proof}

Next, we show that the induced subgraph $P_{2n+1}^k[X_k]$ minimizes the maximum degree among the subgraphs induced by vertex set of size  $\alpha(P_{2n+1}^k)+1$.
\begin{prop}\label{properX2}
Let $H_k=P_{2n+1}^k[X_k]$.        
Then
$$
\Delta(H_k)=
\begin{cases}
	2, &\mbox{ if } n=1,\\
	1, &\mbox{ if } n\geq 2.
\end{cases}
$$
\end{prop}
\begin{proof}
It sufficient to show that $\Delta(H_k) = \Delta(H_{k-1})$. Thus we have 
$$\Delta(H_k)=\Delta(H_{k-1})=\ldots=\Delta(H_1)=\Delta(P_{2n+1}[X_1])=\begin{cases}
	2, &\mbox{ if } n=1,\\
	1, &\mbox{ if } n\geq 2,
\end{cases}$$
as desired.

Recall that
$$
X_k=Q_{2n+1,k}^1(X_{k-1})\cup Q_{2n+1,k}^2(\overline{X_{k-1}})\cup\cdots\cup Q_{2n+1,k}^{2n+1}(X_{k-1}),
$$
and 
$$\overline{X_k}=Q_{2n+1,k}^1(\overline{X_{k-1}})\cup Q_{2n+1,k}^2(X_{k-1})\cup\cdots\cup Q_{2n+1,k}^{2n+1}(\overline{X_{k-1}}).$$
Let $H'_i=P_{2n+1}^{i}[\overline{X_{i}}]$. By Proposition~\ref{neighborcon},
$E(Q_{2n+1,k}^i(X_{k-1}), Q_{2n+1,k}^j(\overline{X_{k-1}}))=\emptyset$ for all $i\in\{1,3,\cdots,2n+1\}$ and $j\in \{2,4,\cdots,2n\}$,
and $E(Q_{2n+1,k}^i(X_{k-1}), Q_{2n+1,k}^j({X_{k-1}}))=\emptyset$ for all $j\not=i$ and $i,j\in\{1,3,\cdots,2n+1\}$. 
Therefore, $H_k$ only has edges in $Q_{2n+1,k}^{i}(X_{k-1})$ and in $Q_{2n+1,k}^{j}(\overline{X_{k-1}})$. 
Hence we have
$$\Delta(H_k)=\max\{\Delta(H_{k-1}),\Delta(H'_{k-1})\}.$$
Similarly, 
we have 
$$\Delta(H'_k)=\max\{\Delta(H_{k-1}),\Delta(H'_{K-1})\}.$$
Thus $\Delta(H_k)=\Delta(H'_k)$ is true for all $k\geq 2$ and hence $\Delta(H_k)=\max\{\Delta(H_{k-1}),\Delta(H'_{K-1})\}=\Delta(H_{k-1})$ for all $k\geq 3$. When $k=1$, it can be checked directly by the definition of $X_1$ that $\Delta(H_1)=1>\Delta(H'_1)=0$. Therefore, $\Delta(H_2)=\max\{\Delta(H_1),\Delta(H'_1)\}=\Delta(H_1)$. 
Consequently, we obtain that
$$
\Delta(H_k)=\Delta(H_{k-1})=\ldots=\Delta(H_1)=
\begin{cases}
	2, &\mbox{if } n=1,\\
	1, &\mbox{if } n\geq 2.
\end{cases}
$$
\end{proof}

Combining Propositions~\ref{properX1} and~\ref{properX2}, we have the following corollary. 
\begin{cor}\label{COR: 2n+1}
For $n\ge 2$,	$$f(P_{2n+1}^k)=1.$$
\end{cor}
By Proposition~\ref{properX2}, we have $f(P_{3}^k)\le2$.
Next, we show that 2 is also the lower bound of  $f(P_3^k)$ using the properties of signed matrices. 
\begin{definition}
Let $$
A_1=\begin{pmatrix}
	0 &1 &0\\
	1 &0 &-1\\
	0 &-1 &0
\end{pmatrix}	
$$
and for $k\ge 1$, define
$$
A_{k+1} = \begin{pmatrix}
	A_k &I &O \\
	I &-A_k &-I \\
	O &-I &A_k
\end{pmatrix}.
$$
\end{definition}

Clearly, 	$A_k$ is a signed matrix of $P_3^k$, with order $3^k$.

\begin{prop}\label{propA2}
The matrix	$A_k$ has eigenvalue $0$ with multiplicity $1$, and the non-zero eigenvalues are pairwise symmetric with respect to $0$. The minimum positive eigenvalue of $A_k$ is $\sqrt{2}.$
\end{prop}
\begin{proof} The proof is by induction on $k$.
When $k=1$, the characteristic polynomial $$P_{A_1}(x)=det\begin{pmatrix}
	x &-1 &O\\
	-1 &x &1\\
	O &1 &x
\end{pmatrix}=x(x+\sqrt{2})(x-\sqrt{2}).$$ The conclusion holds for $k=1$. Now suppose the eigenvalues of $A_k$ are 
$$
\lambda_{3^k}\leq\cdots\leq\lambda_{\frac{3^k+3}{2}}=-\sqrt{2}<\lambda_{\frac{3^k+1}{2}}=0< \sqrt{2}= \lambda_{\frac{3^k-1}{2}}\leq\cdots\leq\lambda_1.
$$
Then the characteristic polynomial of $A_{k+1}$ is 
$$
\begin{aligned}
	P_{A_{k+1}}(x)
	&=\det\begin{pmatrix}
		xI-A_k &-I &O\\
		-I &xI+A_k &I\\
		O &I &xI-A_k
	\end{pmatrix}
	=\det\begin{pmatrix}
		xI-A_k &-I &xI-A_k \\
		-I &xI+A_k &O\\
		O &I &xI-A_k
	\end{pmatrix}\\
	&=\det\begin{pmatrix}
		xI-A_k &-2I &O\\
		-I &xI+A_k &O\\
		O &I &xI-A_k
	\end{pmatrix}
	=\det(xI-A_k)\det\begin{pmatrix}
		xI-A_k &-2I \\
		-I &xI+A_k
	\end{pmatrix}\\
	&=\det(xI-A_k)\det\left[(xI-A_k)(xI+A_k)-2I\right]
	=\det(xI-A_k)\det\left[(x^2-2)I-A_k^2\right]\\
	&=\prod_{i=1}^{3^k}(x-\lambda_{i})\prod_{i=1}^{3^k}(x^2-2-\lambda_{i}^2)
	=\prod_{i=1}^{3^k}(x-\lambda_i)\left(x+\sqrt{2+\lambda_{i}^2}\right)\left(x-\sqrt{2+\lambda_{i}^2}\right).
\end{aligned}
$$
Thus the eigenvalues of $A_{k+1}$ are 
$$\{\lambda_i, \sqrt{2+\lambda_i^2},-\sqrt{2+\lambda_i^2}:i=1,\cdots,3^k\}.$$
By the induction hypothesis of $\lambda_i$, we know that the $A_{k+1}$ has eigenvalue $0$ with multiplicity $1$, the non-zero eigenvalues are pairwise symmetric with respect to $0$ and the minimum positive eigenvalue is $\lambda_{\frac{3^{k+1}-1}{2}}=\sqrt{2}$. 
\end{proof}

\begin{cor}\label{lowP3}
$f(P_3^k)\geq 2$.	
\end{cor}
\begin{proof}
By Theorem~\ref{OddIn}, $\alpha(P_3^k)=\frac{3^k+1}{2}$. Choose arbitrarily an induced subgraph, say $H$, of $P_3^k$, with $\frac{3^k+3}{2}$ vertices. Suppose the principle submatrix of $H$ in $A_k$ is $B$. Then by Lemmas~\ref{Huang1} and~\ref{Cauchy Interlace Thm}
$$
\Delta(H)\geq \lambda_{1}(B)\geq \lambda_{3^k-\frac{3^k+3}{2}+1}(A_k)=\lambda_{\frac{3^k-1}{2}}(A_k)=\sqrt{2}.
$$
Therefore, $f(P_3^k)\geq \lceil\sqrt{2}\rceil=2$.
\end{proof}


\section{Proof of Theorem~\ref{main1} when $m=2n$}

 


We begin our proof by defining the signed matrices of $P_{2n}^k$.
\begin{definition}
Let$$
A_1=\begin{pmatrix}
	0 &1 &0 &\cdots&\cdots &\cdots &0\\
	1 &0 &-1 &\ddots & & &0\\
	0 &-1 &0 &\ddots &\ddots & &0\\
	\vdots & &\ddots &\ddots &\ddots &\ddots&\vdots\\
	0 &\cdots &\cdots &1 &0 &-1 &0\\
	0 &\cdots &\cdots &\cdots &-1 &0 &1\\
	0 &\cdots &\cdots &\cdots &0 &1 &0
\end{pmatrix},
$$	
i.e.,
$$(A_1)_{ij}=\begin{cases}
	1, &\mbox{if } i=2k-1\mbox{, }j=2k\mbox{ or }i=2k\mbox{, }j=2k-1\mbox{, }k\in [n],\\
	-1, &\mbox{if } i=2k\mbox{, }j=2k+1\mbox{ or }i=2k+1\mbox{, }j=2k\mbox{, }k\in [n-1],\\
	0, & \mbox{else}.
\end{cases}
$$
For $k\ge 1$, define 
$$
A_{k+1}=\begin{pmatrix}
	A_k &I \\
	I &-A_k &-I \\
	  &-I &A_k &I \\
	  & &\ddots &\ddots &\ddots \\
	  & & & I &-A_k &-I\\
	  & & & & -I &A_k &I\\
	  & & & & & I &-A_k 
\end{pmatrix}.
$$
\end{definition}

Clearly, $A_{k}$ is a block-tridiagonal signed graph of $P_{2n}^k$. Let $I(k)$ be  a unit matrix with odder $(2n)^k$. It can be checked directly that 	$A_{k}^2=I(1)\otimes A_{k-1}^2+A_1^2\otimes I(k-1)$.



\begin{prop}\label{PropEA4}
	The eigenvalues of $A_k$ are pairwise symmetric with respect to $0$ and $0$ is not an eigenvalue of $A_k$.
\end{prop}
\begin{proof}
	According to Proposition~\ref{CPGBi}, $P_{2n}^k$ is a bipartite graph. Then by Lemma~\ref{bipartite}, all of its non-zero eigenvalues are pairwise symmetric with respect to $0$. 
	We prove by induction on $k$ that $A_k^2>0$ and thus $0$ is not an eigenvalue of $A_k$. 	

	When $k=1$, we show that $rank(A_1)=2n$ and thus $0$ is not an eigenvalue of $A_1$ and, therefore, $A_1^2>0$.
	Let $a_1,a_2\ldots,a_{2n}$ be the row vectors of $A_1$ and $k_1,\ldots,k_{2n}\in \mathbb R$. 
	Then we have  $$
	\begin{aligned}
		\sum_{i=1}^{2n}k_ia_i=0 &\Leftrightarrow
		(k_2,k_1-k_3,-k_2+k_4,\cdots,-k_{2n-1}+k_{2n},k_{2n-1})=0\\
		&\Leftrightarrow k_1=k_3=\cdots=k_{2n-1}=0 
		=k_2=k_4=\cdots=k_{2n}.
	\end{aligned}	
	$$
	Thus $a_1, a_2,\ldots,a_{2n}$ are linearly independent and thus
	$rank(A_1)=2n.$
	
For $k> 1$,	
since $A_1^2>0$ and $A_{k}^2=I(1)\otimes A_{k-1}^2+A_1^2\otimes I(k-1)$, we have $A_{k}^2>0$. 
\end{proof}

If we can find the minimum eigenvalue of $A_1^2$, we can give the minimum eigenvalue of $A_k^2$ by Lemma~\ref{product}. In order to do this, we first give the characteristic polynomial of $A_1^2$.
To begin with, we define some polynomials below:
\begin{definition}
	$$
	f_0(x)=1\mbox{, }f_1(x)=x-2;\mbox{ }f_k(x)=(x-2)f_{k-1}(x)-f_{k-2}(x);
	$$	
	$$
	g_0(x)=1\mbox{, }g_1(x)=x-1;\mbox{ }g_k(x)=(x-2)g_{k-1}(x)-g_{k-2}(x).
	$$
\end{definition}

The following proposition can be checked from the definitions.

\begin{prop}\label{propfg1}
(i)	$x-2-\frac{g_{k-2}}{g_{k-1}}=\frac{g_k}{g_{k-1}}\mbox{; }x-2-\frac{f_{k-2}}{f_{k-1}}=\frac{f_k}{f_{k-1}}.$

(ii) 	$g_{k} = f_k+f_{k-1}$, $k\geq1$.
\end{prop}

Now we calculate the characteristic polynomial of $A_1^2$.
\begin{prop}\label{PropEA1}
$P_{A_1^2}(x)=g_n^2(x)$.
\end{prop}
\begin{proof}
Let  $\Delta_{2n-2i}(f_i(x), g_i(x))=\det\begin{pmatrix}
	\frac{g_{i+1}(x)}{g_{i}(x)} &0 &1 \\
	0 & \frac{f_{i+1}(x)}{f_{i}(x)} &0  &\ddots\\
	1 &0 &x-2 &\ddots &\ddots\\
	&\ddots &\ddots &\ddots &\ddots &\ddots \\
	& &\ddots &\ddots &x-2 &0 &1\\
	& & &\ddots &0 &x-2 &0  \\
	& & & &1 &0 &x-1
\end{pmatrix}.
$ 
Then 

$$
\begin{aligned}
	\det(xI-A_1^2)
	&=\Delta_{2n}(f_0(x), g_0(x))=\det\begin{pmatrix}
		x-1 &0 &1 \\
		0 &x-2 &0  &\ddots\\
		1 &0 &x-2 &\ddots &\ddots\\
		&\ddots &\ddots &\ddots &\ddots &\ddots \\
		& &\ddots &\ddots &x-2 &0 &1\\
		& & &\ddots &0 &x-2 &0  \\
		& & & &1 &0 &x-1
	\end{pmatrix}\\
	&=(x-1)(x-2)\Delta_{2n-2}(f_1(x), g_1(x))=g_1(x)f_1(x)\Delta_{2n-2}(f_1(x), g_1(x))\\
	&=\ldots=g_{n-2}(x)f_{n-2}(x)\Delta_4(f_{n-2}(x),g_{n-2}(x))\\
	&=g_{n-2}(x)f_{n-2}(x)\det\begin{pmatrix}
		\frac{g_{n-1}(x)}{g_{n-2}} &0 &1 &0\\
		0 &\frac{f_{n-1}(x)}{f_{n-2}(x)} &0 &1\\
		1 &0 &x-2 &0\\
		0 &1 &0 &x-1
	\end{pmatrix}\\
	&=g_n(x)f_{n-1}(x)\left(\frac{f_n(x)}{f_{n-1}(x)}+1\right)\\
	&=g_n(x)(f_n(x)+f_{n-1}(x))=g_n^2(x).
\end{aligned}
$$
\end{proof}

\begin{definition}
	Let $\beta_n$ be the minimum positive root of $g_{n}(x)$.
\end{definition}

\begin{prop}\label{PropEA5}
The minimum positive eigenvalue of $A_k$ is $$\lambda_{\frac{1}{2}(2n)^{k}}(A_k)=\sqrt{k\beta_n}.$$
\end{prop}
\begin{proof}
According to Proposition~\ref{PropEA1}, $\beta_n$ is the minimum positive eigenvalue of $A_1^2$. By Proposition~\ref{PropEA4}, $A_1^2>0$. Thus $\beta_n$ is the minimum eigenvalue of $A_1^2$. Note that $A_{k}^2=I(1)\otimes A_{k-1}^2+A_1^2\otimes I(k-1)$. By Lemma~\ref{product} and an inductive argument on $k$, we have   $k\beta_n$ is the minimum eigenvalue of $A_k^2$. Since $A_k$ has no eigenvalue $0$ and the eigenvalues of $A_k$ are pairwise symmetric with respect to $0$, the $\frac{1}{2}(2n)^{k}$-th eigenvalue is the minimum positive one of $A_k$. Therefore, we obtain that $\lambda_{\frac{1}{2}(2n)^{k}}(A_k)=\sqrt{k\beta_n}$.
\end{proof}

Finally, we obtain our main result:
\begin{cor}
	$f(P_{2n}^k)\geq \lceil\sqrt{k\beta_n}\rceil$.
\end{cor}
\begin{proof}
Choose arbitrarily a vertex set $V$ with $\frac{1}{2}(2n)^{k}+1$ vertices and let $H=P_{2n}^k[V].$ Suppose the signed graph of $H$ corresponding  to $A_k$ is $B$. By Lemmas~\ref{Huang1} and~\ref{Cauchy Interlace Thm}, we have 
$$\Delta(H)\geq \lambda_1(B)\geq \lambda_{\frac{1}{2}(2n)^{k}}(A_k)=\sqrt{k\beta_n}.$$
Thus $f(P_{2n}^k)\geq \lceil\sqrt{k\beta_n}\rceil$, where $\beta_n$, a positive constant depending on $n$, is the minimum positive root of $g_n(x)$. 
\end{proof}


\section{Discussion and Remarks}
In this note, we determine the exact value of $f(P_{2n+1}^k)$ and  provide a lower bound of $f(P_{2n}^k)$, but we are not sure whether the lower bound $f(P_{2n}^k)\geq \lceil \sqrt{k\beta_n}\rceil$ is tight for $n\geq 3$. We leave this as a problem.
Another interesting problem proposed by Huang (oral communication):  
\begin{quest}\label{QUES: q2}
	Given a graph $G$, let $f_k(G)$ be minimum of the maximum degree of the induced subgraph with $\alpha(G)+k$ vertices,   what can we say about $f_k(G)$? 
\end{quest}

\section{Acknowledgment}
The work was supported by the National Natural Science Foundation of China (No. 12071453), the National Key R and D Program of China (2020YFA0713100), and the Innovation Program for Quantum Science and Technology, China (2021ZD0302902).


\begin{thebibliography}{99}
\bibitem{Alon}
N. Alon and K. Zheng, Unitary signings and induced subgraphs of Cayley graphs of $Z_n$, Advances in Combinatorics, 2020: 11, 12pp.
	
\bibitem{Chung}
F. R. K. Chung, Z. F\"uredi, R. L. Graham, and P. Seymour, On induced subgraphs of the cube, J. Comb. Theory Ser. A 49 (1988), no. 1, 180-187.

\bibitem{Dong}
D. Dong, On induced subgraphs of the Hamming graph, J. Graph Theory 96 (2021), no. 1, 160-166.

\bibitem{Cauchy}
S. Fisk, A very short proof of Cauchys interlace theorem for eigenvalues of Hermitian matrices, Amer. Math. Monthly, 112 (2005), no.2, 118.

\bibitem{Eigenvalues}
K. A. Germina, S. Hameed K and T. Zaslavsky, On products and line graphs of signed graphs: their
eigenvalues and energy, Linear Alg. Appl., 435 (2011), no.10, 2432-2450. 

\bibitem{Bipartite}
R. Hammack, W. Imrich and S. Klavˇzar, Handbook of product graphs, CRC press, 2011.


\bibitem{AboutConjectur}
P. Hatami, R. Kulkarni, and D. Pankratov, Variations on the Sensitivity Conjecture, Theory Comput. 4 (2011), 1-27.

\bibitem{CartesianType}
Z.-M. Hong, H.-J. Lai, J. Liu, Induced subgraphs of product graphs and a generalization of Huang's theorem,  J. Graph Theory, 98 (2021), no.2, 285-308.

\bibitem{Huang}
H. Huang, Induced subgraphs of hypercubes and a proof of the Sensitivity Conjecture, Ann. Math. 190 (2019), no. 3, 949-955.



\bibitem{S.Sinha}
R. Karthikeyan, S. Sinha, and V. Patil, On the resolution of the sensitivity conjecture, Bull. Amer. Math. Soc. 57 (2020), no. 4, 615-638.

\bibitem{Verret}
F. Lehner, G. Verret, Counterexamples to ``A conjecture on induced sub-graphs of Cayley graphs'', Ars Math. Contemp., 19 (2020), no.1, 77–82.

\bibitem{Tsang}
A. Potechin and H. Y. Tsang, A conjecture on induced subgraphs of Cayley graphs, arXiv:2003.13166.

\bibitem{Vincent}
V. Tandya, An induced subgraph of the Hamming graph with maximum degree 1, J. Graph Theory, 101 (2022), no.2, 311-317.







\bibitem{Tikara}
A. Tikaradze, Induced subgraphs of powers of oriented cycles, Linear Multilinear Algebra, 70 (2020), no.17, 3301-3303.




\end{thebibliography}
\end{document}